\documentclass[a4paper,12pt]{article}
\usepackage{amsthm}
\usepackage{amsmath}
\usepackage{amsfonts}
\usepackage{amssymb}
\usepackage{latexsym}
\usepackage{graphicx}

\usepackage{tabularx}
\usepackage[margin=1.2in]{geometry}
\usepackage{enumitem}
\usepackage{ae}
\usepackage[T1]{fontenc}

\newtheorem{teo}{Theorem}
\newtheorem{prop}{Proposition}

\newtheorem{lem}{Lemma}
\newtheorem{gam}{Game}

\newtheorem{oss}{Remark}
\newtheorem{defin}{Definition}

\title{On MPS construction of blocking sets in projective spaces: a generalization}
\author{ Simone Costa}
\date{}
\begin{document}
\maketitle\begin{abstract}
In this paper we provide a generalization of the MPS construction of blocking sets of
$PG(r,q^n)$ using subspaces of dimension $s\leq n-2$. By this construction, we determine a new non-planar example in $PG(3,q^6)$. 
\end{abstract}
\section{Introduction}
The notion of \emph{blocking set} was introduced for the first time by Richardson in a game theory setting, see \cite{12},
recalling the work of Shapley \cite{13} and of Morgenstern and von Neumann
\cite{11} as a \emph{blocking coalition},
that is a set of players which prevents every coalition
from winning but it is not itself a winning coalition.
Richardson also pointed out
an interesting relationship between the theory of blocking sets and
projective geometry, as made clear by the following example.
\begin{gam}
Consider $PG(2,p^n)$ as a game where the point set $X$ of $PG(2,p^n)$
is the set of players and,
denote by $\mathfrak{F}$ the set of lines of $PG(2,p^n)$.
The minimal winning coalitions are the elements of $\mathfrak{F}$.
In this scenario, the blocking sets turn out to be exactly those subsets
of $X$ which intersect all of the lines, without containing any.
\end{gam}
More generally now, in a projective or affine space, it is usually defined as \emph{k-blocking set} a point set meeting every $(n-k)$-dimensional subspace, and as \emph{blocking set} a point set meeting every hyperplanes.
More generally, in a projective or affine space, it 
a \emph{blocking set with respect to $k$-dimensional subspaces}
is a point set meeting every $k$-dimensional subspace at least once.
A k-blocking set is \emph{trivial} if it contains a $(n-k)$-dimensional
subspace.
In general, any set containing a blocking set is still a blocking set.
Thus, we are interested in \emph{minimal} ones with respect to inclusion.
\par
There is a vast literature on blocking sets (more than $400$ papers on
{\tt mathscinet}). Here, we refer the reader just to \cite{3,4,6,9,10,11}
which are closely related in techniques and topics to the present work.
In \cite{7,8}, Mazzocca, Polverino and Storme have introduced several
constructions of blocking sets of $\Pi_r=PG(r,q^n)$ making use
of the so called
\emph{Barlotti-Cofman representation of $PG(r,q^n)$ in $PG(nr,q)$}.
One of these constructions, which we shall denote as
\emph{MPS construction}, is recalled in Section \ref{sec2}.
In Section \ref{sec3}, a generalization of this
construction is introduced and an example of minimal
blocking set, which appears not to be equivalent to an MPS one, is obtained and
studied in Section \ref{sec6}.
\section{The MPS construction}
\label{sec2}
By MPS construction we mean the construction carried out by Mazzocca,\ Polverino, Storme in \cite{7}: starting from a blocking set in a projective space, one can construct blocking sets in spaces whose order is a power of the original one. The idea of the construction generalizes the planar version of Mazzocca and
Polverino in \cite{8}; in this section we follow \cite{5}.\\
We consider the \emph{Barlotti-Coffman} representation
of $PG(r,q_1^n)$ in $PG(nr,q_1)$.
Take $\Sigma=PG(nr-1,q_1)$ and let $\mathcal{S}$ be one of its desarguesian
$(n-1)$-spreads; see \cite{1,5bis} and also \cite{2}.
\begin{defin}
Embed $\Sigma$ as an hyperplane in $\Sigma'=PG(nr,q_1)$ and define a \textbf{point-line geometry} $\mathbf{\Pi_r=\Pi_r(\Sigma',\Sigma,\mathcal{S})}$ in the following way:
\begin{itemize}
\item the points of $\Pi_r$ are the points of $\Sigma'\setminus \Sigma$ and the elements of $\mathcal{S}$; 
\item the lines of $\Pi_r$ are the n-subspaces of $\Sigma'$ intersecting $\Sigma$ in an element of $\mathcal{S}$, and the lines of $PG(\mathcal{S})$;
\item the point-line incidences are inherited from $\Sigma$ and $\Sigma'$.
\end{itemize}
\end{defin}
\begin{teo}[\cite{1Bis}]
The incidence structure $\Pi_r$ is isomorphic to the projective space $PG(r,q)$, where $q=q_1^n$.
\end{teo}
This incident structure is called
the {\bfseries{Barlotti-Cofman representation of $PG(r,q)$}}.\\
The points of $\Pi_r$ in $\Sigma'\setminus \Sigma$ are called \emph{affine}.
Let $Y$ be a fixed element of $\mathcal{S}$ and let $\Omega=\Omega_{n-2}$ be an hyperplane of $Y$. Let $\Gamma'=\Gamma'_{(r-1)n+1}$ be an $((r-1)n+1)$-subspace of $\Sigma'$ disjoint from $\Omega$. Also, denote by $\Gamma=\Gamma_{(r-1)n}$ the $(r-1)n$-subspace intersection of $\Gamma'$ and $\Sigma$ and by $t$ the intersection point of $\Gamma$ and $Y$.\\
Let $\overline{B}$ be a blocking set of $\Gamma'$ such that $\overline{B}\cap \Gamma=\{\tilde{q}\}$, $\tilde{q}$ a point, with the following property:
$$L\setminus \{t\}\not \subset \overline{B},$$
for every line $L$ of $\Gamma'$ through $t$.\\
Denote by $K=\mathfrak{K}(\Omega,\overline{B})$ the cone with vertex $\Omega$ and base $\overline{B}$. Note that, since $\Gamma' \cap \Omega=\emptyset$, \ we have \ $\left\langle \overline{p},\Omega\right\rangle\cap \left\langle \overline{p'},\Omega\right\rangle=\Omega$, for any two distinct points $\overline{p'},\overline{p}\in \overline{B}$. Let $B$ be the subset of $\Pi_r$ defined by 
$$ B=(K\setminus \Sigma)\cup \{ X \in \mathcal{S}: X \cap K \not= \emptyset\}, $$
and note that
\begin{itemize}
\item if $\tilde{q} \in Y$ (i.e. $t=\tilde{q}$), then $|B|=q_1^{n-1}(|\overline{B}|-1)+1$ and $B\cap PG(\mathcal{S})=\{Y\}$;
\item while if $\tilde{q}\not\in Y $, then $|B|=q_1^{n-1}|\overline{B}|+1$ and $|B\cap PG(\mathcal{S})|=q_1^{n-1}+1$.
\end{itemize}
In both cases we get:
\begin{teo}[\cite{7}] 
$B$ is a blocking set of the projective space $\Pi_r$.
\end{teo}
Suppose that $\overline{B}$ is a minimal blocking set of $\Gamma'$ such that $t=\tilde{q}$ (this case is called MPS construction A), in other words, suppose that $\Gamma$ is a tangent hyperplane of $\overline{B}$ at the point $\tilde{q}$.\\
In this case, $B=(K\setminus \Sigma)\cup \{Y\}$ and $$ |B|=|K\setminus\{\Sigma\}|+1= q_{1}^{n-1}(|\overline{B}|-1)+1. $$
Then the following theorem holds true; see \cite{7}.
\begin{teo}
$B$ is a minimal blocking set of $\Pi_r$ if and only if $\overline{B}$ is a minimal blocking set of $\Gamma'$.
\end{teo}
\section{The construction}
\label{sec3}
Let us consider $\Pi_r=PG(r,q^n)$ as represented (using the Barlotti-Cofman representation) into $PG(nr,q)$ with respect to a desarguesian spread $\mathcal{S}$ of an hyperplane $\Sigma$ of $PG(nr,q)$.\\
Let $X$ be a fixed element of $\mathcal{S}$ and $\Omega$ be a subspace of $X$ of dimension $s\leq n-2$.\\
Let $\Gamma$ be a subspace of $\Sigma$ such that $\ \dim \Gamma= rn-s-2$ and \ $\Gamma \cap \Omega =\emptyset$, let $\Theta= \Gamma \cap X$ and $\Gamma'$ be a subspace of $PG(nr,q)$ of dimension $rn-s-1$ such that $\Gamma'\cap \Sigma=\Gamma$.\\
Let $\mathfrak{F}=\mathfrak{F}_{(r-1)n}^s=\{\left\langle S_{(r-1)n},\Omega\right\rangle\cap \Gamma'$ : $S_{(r-1)n}$ hyperplane of $\Pi_r$ not containing $X$\}.
Because $\Gamma'\cap \Omega=\emptyset$ by the Grassmann formula we get:
$$\dim\left(\left\langle \Gamma',\Omega\right\rangle\right)-1=rn-s-1+s $$ hence
$$\dim\left(\left\langle \Gamma',\Omega\right\rangle\right)=rn= \dim\left(\left\langle \left\langle \Gamma',\Omega\right\rangle, S_{(r-1)n}\right\rangle\right) $$
and, called $I_{(r-1)n}$ an element of $\mathfrak{F}$, we get:
$$\dim\left(\left\langle \left\langle \Omega,S_{(r-1)n}\right\rangle, \Gamma'\right\rangle\right)+\dim\left(\left\langle \Omega,S_{(r-1)n}\right\rangle\cap \Gamma'\right)=rn+\dim(I_{(r-1)n})= $$
$$=\dim\left(\left\langle \Omega,S_{(r-1)n}\right\rangle\right)+\dim(\Gamma')=(r-1)n+s+1+rn-s-1=2rn-n.$$
It follows that $\dim(I_{(r-1)n})=(r-1)n$ and $\mathfrak{F}$ is a family of subspaces of $PG(rn,q)$ of dimension $(r-1)n$.
\begin{defin}
A subset $\overline{B}$ of $\Gamma'$ is called \textbf{blocking set with respect to $\mathfrak{F}$} or simply an $\mathfrak{F}$-\textbf{blocking set} if given $ I_{(r-1)n}\in \mathfrak{F}$, we get $\overline{B}\cap I_{(r-1)n} \not = \emptyset.$ An $\mathfrak{F}$-\textbf{blocking set} is called minimal if it is minimal with respect to the inclusion.\\
Moreover, in the case $\mathfrak{F}$ be a family of subspaces of $PG(m,q)$ of dimension $k$ and $B$ an $\mathfrak{F}$-blocking set, we call $B$ \emph{trivial} if it contains a subspace of dimension $m-k$.
\end{defin}
Let $\overline{B}$ be an $\mathfrak{F}$- blocking set of $\Gamma'$ such that $\overline{B}\cap \Sigma= \Theta$, then we consider the cone $K=\mathfrak{K}(\Omega, \overline{B})$ and $B=K\cup \{X\}$. We call this construction ``\textbf{generalized MPS construction}''.
\begin{teo}
$B\cap S_{(r-1)n}\not= \emptyset$ for all hyperplanes $S_{(r-1)n}$ of $\Pi_r$.
\end{teo}
\proof
Consider an hyperplane $S_{(r-1)n}$ not containing $X$. Hence there exists $p' \in\left\langle S_{(r-1)n},\Omega\right\rangle\cap \overline{B}$ and $p \in \mathfrak{K}(\Omega, \overline{B}) \cap S_{(r-1)n}$.
\endproof
Now we introduce some useful lemmas for proving the minimality of $B$.
\begin{lem}
Let $P_1,P_2$ be subspaces of a projective space $\Gamma$ such that $P_1\cap P_2=\emptyset$.\\
Then $\left\langle \overline{p},P_1\right\rangle\cap\left\langle \overline{p}',P_1\right\rangle=P_1$ for any two distinct points $\overline{p},\ \overline{p}'\in P_2$.
\end{lem}
\proof
Suppose there is an intersection point outside $P_1$, this means that $\langle \overline{p},P_1 \rangle=\langle\overline{p'},P_1 \rangle $. Then the line through $\overline{p}$ and $\overline{p'}$ meets $P_1$, which implies that $P_2$ meets $P_1$ non-trivially, a contradiction.
\endproof
\begin{lem}
Let us consider $\overline{B}\setminus X$, then for $I_{(r-1)n}=\left\langle S_{(r-1)n},\Omega\right\rangle\cap \Gamma' \in \mathfrak{F}$ and $S_{(r-1)n}$ not containing $X$ we get:
$$ |B \cap S_{(r-1)n}|=|(\overline{B}\setminus{X}) \cap I_{(r-1)n}|. $$ 
\end{lem}
\proof
Consider:
$$B \cap S_{(r-1)n}=(K\setminus \Sigma) \cap S_{(r-1)n}=\bigcup_{\overline{p}\in \overline{B}\setminus{X}}(\left\langle \overline{p},\Omega\right\rangle\cap S_{(r-1)n})$$
Since $\ \dim\left(\left\langle \Omega,S_{(r-1)n}\right\rangle\right)=(s+1)+(r-1)n\ $ and, for $\overline{p}\in \overline{B}\setminus{X}$
$$\dim\left(\left\langle \overline{p},\Omega\right\rangle\cap S_{(r-1)n}\right)+\dim\left(\left\langle \overline{p},\Omega,S_{(r-1)n}\right\rangle\right)=(s+1)+(r-1)n,$$
we get $\dim\left(\left\langle \overline{p},\Omega\right\rangle\cap S_{(r-1)n}\right)\in \{-1,0\}$ and equals $0$ if and only if $\overline{p}\in \langle \Omega,S_{(r-1)n} \rangle$ and since $\overline{p}\in \overline{B}\setminus{X}\subseteq \Gamma'\setminus \Sigma$ this holds if and only if $\overline{p}\in(\overline{B}\setminus{X})\cap I_{(r-1)n}$.\\
Since (by Lemma 1) $\left\langle \overline{p},\Omega\right\rangle\cap\left\langle \overline{p}',\Omega\right\rangle=\Omega$ for any two distinct points $\overline{p},\ \overline{p}'\in \overline{B}$ we get:
$$\bigcup_{\overline{p}\in \overline{B}\setminus{X}}(\left\langle \overline{p},\Omega\right\rangle\cap S_{(r-1)n})=\bigsqcup_{\overline{p}\in (\overline{B}\setminus{X})\cap I_{(r-1)n}}(\left\langle \overline{p},\Omega\right\rangle\cap S_{(r-1)n}).$$
The claim follows by a simple counting.
\endproof
\begin{teo}
Let $\overline{B}\setminus X$ be a minimal (and such that $\overline{B}$ is non-trivial) $\mathfrak{F}$-blocking set, then $B$ is a minimal (non-trivial) blocking set of $\Pi_r$.
\end{teo}
\proof
Let $\overline{B}\setminus X$ be a minimal blocking set, then by Theorem 4, $(B\setminus \{X\})\cap S_{(r-1)n}'\not = \emptyset$ for all hyperplanes $S_{(r-1)n}'$ not containing $\{X\}$.\\
Consider now a point $p\in B\setminus \{X\}$ then there exists at least one point $\tilde{q}$ such that:
$$\tilde{q}\in \overline{B}\ \cap \left\langle p,\Omega\right\rangle.$$
Let now $\tilde{q},\tilde{q}'$ be such that $\tilde{q}, \tilde{q}'\in \overline{B}\ \cap \left\langle p,\Omega\right\rangle$ which implies $\left\langle p,\Omega\right\rangle=\left\langle \tilde{q},\Omega\right\rangle=\left\langle \tilde{q}',\Omega\right\rangle.$ Since, for all $\tilde{q}\not= \tilde{q}'\in \overline{B}$ we have $ \left\langle \tilde{q},\Omega\right\rangle\cap \left\langle \tilde{q}',\Omega\right\rangle=\Omega$ we get $\tilde{q}=\tilde{q}'$. Hence the point $\tilde{q}$ is uniquely determined by the point $p$.\\
Because of the minimality of $\overline{B}\setminus X$ there exists $I_{(r-1)n}=\left\langle \Omega,\overline{S}_{(r-1)n}\right\rangle\cap\ \Gamma'$
such that $\ \{\tilde{q}\}= I_{(r-1)n}\cap \overline{B}\setminus X.$ We consider now the hyperplane of $\Pi_r$ represented by: $$S_{(r-1)n}=\left\langle p, \overline{S}_{(r-1)n}\cap \Sigma\right\rangle.$$ 
Then, because $\tilde{q}\in\langle p,\Omega\rangle$ and $\tilde{q}\in\langle \Omega, S_{(r-1)n}\rangle$ we have: 
$$\left\langle \Omega,S_{(r-1)n}\right\rangle\cap\ \Gamma'=\left\langle \Omega,p,\overline{S}_{(r-1)n}\cap \Sigma\right\rangle\cap\ \Gamma'=$$
$$=\left\langle \Omega,\tilde{q},\overline{S}_{(r-1)n}\cap \Sigma\right\rangle\cap\ \Gamma'=\left\langle \Omega,\overline{S}_{(r-1)n}\right\rangle\cap\ \Gamma'=I_{(r-1)n}.$$
Hence, since $|I_{(r-1)n}\cap (\overline{B}\setminus X)|=|\{\tilde{q}\}|=1$, because of Lemma 2, $S_{(r-1)n}$ intersects $B$ exactly in the point $p$.\\
Suppose now $B$ contains a line of $\Pi_r$, i.e.\ $\exists \ S_n: S_n\subseteq B$. Hence $S_n=B $ by minimality and $S_n= \bigcup_{\overline{p}\in \overline{B}} \left\langle \overline{p},\Omega\right\rangle$. 
Since $\Gamma'\cap \Omega=\emptyset$ for $\overline{p}\in \overline{B}\subseteq \Gamma'$ we have $\Gamma'\cap \langle \overline{p},\Omega\rangle =\overline{p}$ and hence:
$$\Gamma'\cap S_n=\Gamma'\cap (\cup_{\overline{p}\in \overline{B}}\langle \overline{p},\Omega\rangle)=\cup_{\overline{p}\in \overline{B}}(\Gamma'\cap\langle \overline{p},\Omega\rangle)=\cup_{\overline{p}\in \overline{B}}\overline{p}=\overline{B}. $$
Therefore if $B$ is trivial then $\overline{B}$ is trivial.
\endproof
\begin{teo}
In the situation of the previous theorem the following equality holds:
$$|B|=(|\overline{B}|-\frac{q^{n-s-1}-1}{q-1})q^{s+1}+1 $$
\end{teo}
\proof
For each $\overline{p},\overline{p'}\in \overline{B}\setminus \{X\}$ we have:
$$\left\langle \overline{p},\Omega\right\rangle\cap\left\langle \overline{p}',\Omega\right\rangle=\Omega;$$
because $B\cap PG(\mathcal{S})=\{X\}$ we get:
$$|B|=(|\overline{B}|-|\overline{B}\cap \Sigma|)(|\left\langle \overline{p},\Omega\right\rangle|-|\Omega|)+1$$
and because $\overline{B}\cap \Sigma=\Theta$ the claim follows easily.
\endproof
\begin{oss}
If $r=2$ the above construction coincides with the second construction given by Mazzocca and Polverino in \cite{8} (MP Construction $B$); if $s=n-2$ the above construction coincides with the first construction given by Mazzocca, Polverino and Storme in \cite{7} (MPS Construction $A$).
\end{oss}
Our goal is now to find some example of minimal $\mathfrak{F}$-blocking sets $\overline{B}\setminus X$ in order to use our construction and obtain minimal blocking sets with respect to hyperplanes.
\section{Non-planar example}
\label{sec6}
Let consider $\Pi_3=PG(3,q^{6})$ as represented (via Barlotti-Cofman representation) with respect to the desarguesian 2-spread $\mathcal{S}$ of an hyperplane $\Sigma$ of $PG(9,q^2)$.
Let $X,X'$ be two fixed elements of $\mathcal{S}$ and let $p$ be a point of $X$.
Let $\Gamma$ be a subspace of $\Sigma$ of dimension $7$ such that $X'\subseteq \Gamma$ and $p\not \in \Gamma$. Let $\Theta=\Gamma\cap X=\langle r,\tilde{q}\rangle$ for $r,\tilde{q} \in \Theta$ and let $\Gamma'$ be a subspace of $PG(9,q^2)$ of dimension 8 such that $\Gamma'\cap \Sigma=\Gamma$.\\
Let $\mathfrak{F}:= \{\left\langle S_{6},p\right\rangle\cap\ \Gamma':\ S_{6} \mbox{ is an hyperplane of }\Pi_3 \mbox{ not containing } X\}.$
Moreover let us consider the family of seven dimensional subspaces of $PG(9,q^2)$ defined by:
$$\mathfrak{H}:= \{\left\langle S_{6},p\right\rangle:\ S_{6} \mbox{ is an hyperplane of }\Pi_3 \mbox{ not containing } X\}.$$
It is clear that a subset of $\Gamma'$ is a minimal $\mathfrak{F}$-blocking set if and only if it is a minimal $\mathfrak{H}$-blocking set.\\
Let $\pi$ be a plane of $\Gamma'$ through the point $\tilde{q}$ and a point $t\in X'$, not contained in $\Sigma$, that is $\pi\cap \Sigma=\langle t,\tilde{q}\rangle$.
Let $V$ be a Baer subplane of $\pi$ such that $V\cap \langle t,\tilde{q}\rangle =\{\tilde{q}\}$.\\
We recall that every line $L$ of $\pi$ intersects $V$ in either $1$ or $q+1$ points: in the first case we say that $L$ is an \textbf{imaginary} line for $V$ and in the second case we say that $L$ is a \textbf{real} line for $V$. It is possible to prove that for each point $u$ of $\pi\setminus V$ there exists an unique real line $L:\ u\in L.$ 
\\
Name by $S_3$ the 3-space containing $\pi$ and $r$, we consider the Baer cone $\mathfrak{K}(r,V)$ and we observe that every line $L$ of $S_3$ not through the point $r$, is either a real line or an imaginary line for some Baer subplane contained in the cone $\mathfrak{K}(r,V)$ and so intersects the cone in either $1$ or $q+1$ points: in the first case we say that $L$ is an \textbf{imaginary} line for the cone $\mathfrak{K}(r,V)$ and in the second case we say that $L$ is a \textbf{real} line for $\mathfrak{K}(r,V)$.\\
Let us consider a real line $L$ of the plane $\pi$ through the point $t\in X'$, a point $s\in V\cap L$ and we construct the cone:
$$\bar{B}:=(\mathfrak{K}(r,V)\setminus\mathfrak{K}(r,L))\cup\mathfrak{K}(r,s)=\mathfrak{K}(r,V\setminus L\cup \{s\}).$$
Then we have the following Proposition.
\begin{prop}
The cone $\bar{B}$ is a blocking set with respect to the family of subspaces
$$ \mathfrak{H}_t:=\{S_7\in \mathfrak{H}:\ X'\subseteq S_7\}.$$
\end{prop}
\proof
Since the cone $\mathfrak{K}(r,V)$ is a Baer Cone of type $(2,0)$, it is a blocking set with respect to the seven dimensional subspaces of $PG(9,q^2)$ (see \cite{5},\cite{5Tris}) and hence a $\mathfrak{H}$-blocking set. Therefore the elements of $\mathfrak{H}_t$ are in the form $\left \langle p,u, X',Y_1\right \rangle$ with $Y_1\in \mathcal{S},\ Y_1\not\subset \left \langle X',X\right \rangle$ and $u\in \mathfrak{K}(r,V) $. Then, if $u\in \bar{B}$, we get that $\left \langle p,u, X',Y_1\right \rangle\cap \bar{B}\not=\emptyset$. Let us suppose $u\not \in \bar{B}$, hence $u\in \mathfrak{K}(r,L\setminus \{s\})$. But in this case $\langle t,u\rangle\cap \mathfrak{K}(r,s)\not=\emptyset$ and hence, since $t\in X'$:
$$\left \langle p,u, X',Y_1\right \rangle\cap \bar{B}\supseteq \langle t,u\rangle\cap \bar{B}\supseteq \langle t,u\rangle\cap \mathfrak{K}(r,s) \not=\emptyset.$$
We conclude that $\bar{B}$ is an $\mathfrak{H}_t$-blocking set.
\endproof
Now we construct a set $\tilde{B}$ such that $\bar{B}\cup \tilde{B}$ is a non-planar, minimal $\mathfrak{H}$-blocking set. For doing this we need the following lemma.
\begin{lem}
There exists exaclty one point $\tilde{t}\in X',\ \tilde{t}\not =t$ such that for all $S_2\in \mathcal{S}, \ S_2\subset \left \langle X,X' \right \rangle$ with $\tilde{t}\in \langle p,S_2 \rangle$ we have $\left \langle p, S_2 \right \rangle \cap \langle t,r\rangle\not=\emptyset.$
\end{lem}
\proof
Since the spread $\mathcal{S}$ is desarguesian the set $\mathfrak{R}:=\{S_2\in \mathcal{S}:S_2\cap \langle t,r\rangle \not=\emptyset\}$ is a regulus. Clearly $X\in \mathfrak{R}$ and hence there exists exactly one line $l$ such that $p\in l$ and $\mathfrak{R}=\{S_2\in \mathcal{S}: S_2\cap l\not=\emptyset\}$. Since also $X'\in \mathfrak{R}$ we have $X'\cap l\not=\emptyset$; let we call $\tilde{t}=X'\cap l$. Then, $l=\langle p,\tilde{t}\rangle$ and hence:
$$\mathfrak{R}=\{S_2\in \mathcal{S}: S_2\cap l\not=\emptyset\}=\{S_2\in \mathcal{S}: \tilde{t}\in \langle p,S_2 \rangle\} .$$
But, for the definition of $\mathfrak{R}$ this means that for all $S_2\in \mathcal{S}:\tilde{t}\in \langle p,S_2\rangle$ we have:
$$\emptyset\not=S_2\cap \langle t,r\rangle \subseteq \langle p,S_2\rangle\cap \langle t,r\rangle.$$
Let now consider $S_2\in \mathcal{S}$ such that $\langle p,S_2\rangle \cap \langle t,r\rangle\not=\emptyset$ and let $u\in \langle p,S_2\rangle \cap \langle t,r\rangle$. Let us consider $S_2'\in \mathcal{S}$ such that $u\in S_2'$; then, since $\langle p,u\rangle \subseteq\langle p,S_2\rangle \cap \langle p,S_2'\rangle$ we must have $S_2=S_2'$. Therefore $u\in S_2 \cap \langle t,r\rangle$ and hence $S_2\in \mathfrak{R}$ and
$$\mathfrak{R}=\{S_2\in \mathcal{S}: \langle p,S_2\rangle \cap \langle t,r\rangle\not=\emptyset\}.$$
Since for a point $y\in X'$ different from $\tilde{t}$ the regulus 
$$\mathfrak{R'}:=\{S_2\in \mathcal{S}:y\in \langle p,S_2 \rangle\}$$
is different from $\mathfrak{R}$, and hence contains an element not contained in $\mathfrak{R}$, we have that there exists an unique required point $\tilde{t}\in X'.$
\endproof
Let us consider lines $L_1,L_2,L_3\subseteq X'$ such that $L_1\cap L_2\cap L_3=\emptyset$ and $t,\tilde{t}\not \in L_i\ i\in \{1,2,3\}$ where $\tilde{t}$ is the point of the previous Lemma. Let us consider a point $h\in \Gamma',\ h\not \in \left \langle X,X',V\right \rangle$, then we define:
$$\tilde{B}:=\mathfrak{K}(h, L_1\cup L_2\cup L_3)\setminus \Sigma. $$
Then:
\begin{prop}
$\tilde{B}$ is a blocking set with respect to the family of subspaces
$$\mathfrak{H}\setminus \mathfrak{H}_t:=\mathfrak{H}\setminus \{S_7\in \mathfrak{H}:\ X'\subseteq S_7\}.$$
\end{prop}
\proof
Given $S_7 \in \mathfrak{H}\setminus \mathfrak{H}_t,$
we have:
$$S_7 \cap \mathfrak{K}(h,X')\supseteq R\equiv PG(1,q^2):\ R\not \subset X'.$$
Since $\mathfrak{K}(h,X')\setminus X'\equiv A(3,q^2)$ and the union of three planes, two by two non-parallel, is a blocking set with respect to the lines in $A(3,q^2)$ we have that $\tilde{B}$ is a blocking set with respect to $\mathfrak{H}\setminus \mathfrak{H}_t.$
\endproof
\subsection{Minimality}
Now we characterize some property of the spectrum of intersection between $\bar{B}\cup \tilde{B}$ and the elements of $\mathcal{S}$. Let we call $S_3:=\langle \bar{B}\rangle$ and $S_2:=S_3\cap \Sigma$.
\begin{prop}
Let $S_7\in \mathfrak{H}$, then $S_7\cap S_3$ is a line not contained in $\Sigma$ and hence $|S_7\cap \bar{B}|\in \{0,1,q,q+1\}$. 
\end{prop}
In particular if $S_7\cap S_3$ is a real line of $S_3$ (w.r.t. the Baer Cone $\mathfrak{K}(r,V)$) through the point $t$, then it is contained in the plane $\langle t,r,s\rangle$.
\proof
Let $S_7\in \mathfrak{H}$, then $S_7=\langle p,Z,Z',n\rangle$ for some $n\not\in \Sigma$ and for some $Z,Z'\in \mathcal{S}$ with $X\not\subseteq \langle Z,Z' \rangle$.\\
For the Grassmann formula we have that $dim(S_3\cap S_7)\geq 1$. Let us suppose that $dim(S_3\cap S_7)\geq 2$ which means that there exists a plane $\pi_2\subseteq S_3\cap S_7$ and therefore a line $L\subseteq (S_3\cap S_7)\cap \Sigma$; since $p\not\in S_3$ we get $p\not\in L$.\\
We have that $L\subseteq S_3\cap \Sigma=\langle r,\tilde{q},t\rangle$ and hence that $L\cap \langle r,\tilde{q} \rangle\subseteq L\cap X \not=\emptyset$. But since 
$$\{p\}=S_7\cap X\supseteq L\cap X\not= \emptyset$$
and $p\not\in L$ we get a contradiction. Therefore $dim(S_3\cap S_7)=1$ and, similarly, we have $dim((S_3\cap S_7)\cap \Sigma)=0$.\\
Let $L_{S_7}=S_7\cap S_3$ and note that $r\not\in L_{S_7}$ since $X\not\subseteq S_7$.
By definition, if the line $L_{S_7}$ is imaginary then $|S_7\cap \mathfrak{K}(r, V) |$ is 1 and if the line $L_{S_7}$ is real $|S_7\cap \mathfrak{K}(r, V) |=q+1$. 
Then we can get the intersection with $\bar{B}$:
if the line $L_{S_7}$ does not intersect $\mathfrak{K}(r, V)\setminus \mathfrak{K}(r, (V\setminus L)\cup \{s\})$ the intersection is equal to $|S_7\cap \mathfrak{K}(r, V) |$ otherwise if $L_{S_7}\not\subseteq \langle r,L\rangle$ the intersection decreased by one. Lastly, if $L_{S_7}\subseteq \langle r,L\rangle$ the intersection is always one.\\
Therefore we have the following possibilities for the intersection $L_{S_7}\cap \bar{B}$:
\begin{center}
\begin{tabular}{|l|r|}
\hline {$\mbox{ Line } S_7\cap \langle \bar{B}\rangle $} & $S_7\cap \overline{B}$ \\
\hline {$\mbox{ Imaginary }$} & $0$ \\
\hline {$\mbox{ Real }$} & $1$ \\
\hline {$\mbox{ Imaginary }$} &$1$ \\
\hline {$\mbox{ Real }$} & $q+1$ \\
\hline {$\mbox{ Real }$} & $q$ \\
\hline
\end{tabular}
\end{center}
\endproof
\begin{prop}
Let $S_7\in \mathfrak{H}_t$, then:
\begin{itemize}
\item $S_7\cap S_3$ is a line through $t$ not contained in $\Sigma$.
\item $S_7\cap \tilde{B}=\emptyset$ or $S_7\cap \tilde{B}=\tilde{B}$.
\end{itemize}
\end{prop}
\proof
Since $X'\subseteq S_7$, $t\in S_7\cap S_3$, hence by Prop. 3, the intersection $S_7\cap S_3$ is a line through $t$ not contained in $\Sigma$.\\
Let $\tilde{S_3}=\langle \tilde{B}\rangle$, i.e. $\tilde{S_3}=\langle X',h\rangle$.
Then either $\tilde{S}_3\subseteq S_7$ or $\tilde{S}_3\cap S_7=X'$. In the first case $S_7\cap \tilde{B}=\tilde{B}$, in the second case $S_7\cap \tilde{B}=\emptyset$.
\endproof
On the other hand for an element $S_7\in \mathfrak{H}\setminus \mathfrak{H}_t$ (which means that $X'\not \subseteq S_7$) we have:
\begin{prop}Let $S_7\in \mathfrak{H}$ such that $X'\not \subseteq S_7$. Then we have that
$S_7\cap \langle \tilde{B}\rangle$ is a line $L'_{S_7}$ not contained in $\Sigma$ and hence $|S_7\cap \tilde{B}|\in \{1,2,3,q^2\}$.
\end{prop}
\proof
Let $S_7=\langle p,u,Z,Z'\rangle$ for some $Z,Z'\in \mathcal{S}$ such that $X,X'\not \subseteq \langle Z,Z'\rangle$ and $u\not \in \Sigma$. Since $\langle u,Z,Z'\rangle$ can be seen as a plane of $\Pi_3$ and $\langle h,X'\rangle$ as a line not contained in $\langle u,Z,Z'\rangle$, we have that $\langle u,Z,Z'\rangle \cap \langle h,X'\rangle=\langle u,Z,Z'\rangle \cap\langle \tilde{B}\rangle$ is just a point. Therefore $S_7\cap \langle \tilde{B}\rangle$ is a line $L'_{S_7}$ not contained in $\Sigma$. Now, since $\tilde{B}$ is the union of three non parallel affine planes we have $|S_3\cap \tilde{B}|=|L'_{S_7}\cap \tilde{B}|\in \{1,2,3,q^2\}$.
\endproof
Hence if $S_7\in\mathfrak{H}$, by Prop. 4 and 5 we have the following possibilities:
\begin{center}
\begin{tabular}{|l|r|}
\hline {$\mbox{Intersection }$} & $S_7\cap \tilde{B} $ \\
\hline {$\emptyset$} & $0 $ \\
\hline {$\langle \tilde{B}\rangle$} & $|\tilde{B}| $ \\
\hline {$\mbox{ Line }$} & $q^2 $ \\
\hline {$\mbox{ Line }$}& $3 $ \\
\hline {$\mbox{ Line }$} & $2 $ \\
\hline {$\mbox{ Line }$} & $1 $ \\
\hline
\end{tabular}
\end{center}
Now we apply this result.
\begin{prop}
Let $u\in \bar{B}$, then there exists an element $S_7\in \mathfrak{H}_t$ tangent to $\bar{B}\cup \tilde{B}$ in $u$. 
\end{prop}
\proof
Let us suppose $u\in \bar{B}\setminus \{s\}$, we determine an element of the family $\mathfrak{H}$ tangent to $u$. Because in $\Pi_3$ exists a plane through $s,h,\{X'\}$, there exists $Y\in \mathcal{S}$ such that $h\in \left\langle s,Y,X' \right \rangle$. In particular $Y\not \subset \langle X,X'\rangle$. 
Then, called $S_7=\left \langle p,u,Y,X' \right \rangle$ we have $S_7\in \mathfrak{H}_t$ and hence, for Proposition 4:
$$S_7 \cap S_3=\langle t,u\rangle.$$
If the line $\langle t,u\rangle$ is imaginary then we have:
$$S_7 \cap \bar{B}=\{u\}.$$
If the line $\langle t,u\rangle$ is real then we have $\langle t,u\rangle\subseteq \left \langle t,s,r \right \rangle$ and hence:
$$S_7 \cap \bar{B}=(S_7\cap S_3) \cap \bar{B}=$$
$$=\langle t,u\rangle\cap (\bar{B}\cap\left \langle t,s,r \right \rangle) =\langle t,u\rangle\cap \langle r,s\rangle=\{u\}.$$
For proving that $S_7$ is tangent we have to estimate the intersection with $\tilde{B}$.
Because of the choice of $Y$, we have $$\left\langle p,h,Y,X'\right \rangle=\left\langle p,s,Y,X'\right \rangle,$$
and, since we have proven that $ \bar{B}\cap S_7=\{u\}$, we have $s\not\in S_7$ and:
$$\left\langle p,h,Y,X'\right \rangle=\left\langle p,s,Y,X'\right \rangle\not =\left\langle p,u,Y,X'\right \rangle=S_7$$
and hence $h\not \in S_7$. Therefore we have that $S_7\cap \tilde{B}\not=\tilde{B}$ and, being $S_7\in \mathfrak{H}_t$, for Proposition 4 we have $S_7\cap \tilde{B}=\emptyset$.
Hence $S_7$ is tangent to $\bar{B}\cup \tilde{B}$ in $u$.\\
Let now $u=s$, and let $s'\in \bar{B}\setminus \{X'\},\ s'\not=s$.
Because, in $\Pi_3$ exists a plane through $s',h,\{X'\}$, there exists $Y_1\in \mathcal{S}\setminus\{X'\}$ such that $h\in \left\langle s',Y_1,X' \right \rangle$. 
Called, as before, $S_7=\left\langle p,s,Y_1,X'\right \rangle$, we have that $S_7\in \mathfrak{H}_t$ and, for Proposition 4:
$$S_7 \cap \bar{B}=(S_7 \cap S_3)\cap \bar{B}=\langle t,s\rangle\cap \bar{B}=\langle t,s\rangle\cap \bar{B}=\{s\}.$$
For proving that $S_7$ is tangent we have to estimate again the intersection with $\tilde{B}$. 
Because of the choice of $Y_1$, we have $$\left\langle p,h,Y_1,X'\right \rangle=\left\langle p,s',Y,X'\right \rangle$$
and since we have proven that $ \bar{B}\cap S_7=\{s\}$ and hence $s'\not\in S_7$ we have:
$$\left\langle p,h,Y_1,X'\right \rangle=\left\langle p,s',Y_1,X'\right \rangle\not =\left\langle p,s,Y_1,X'\right \rangle=S_7$$
and hence $h\not \in S_7$. Therefore we have that $S_7\cap \tilde{B}\not=\tilde{B}$ and, beeing $S_7\in \mathfrak{H}_t$, for Proposition 4 we have $S_7\cap \tilde{B}=\emptyset$.
Hence we can find a tangent element of $\mathfrak{H}$ to all $u\in\bar{B}$
\endproof
\begin{prop}
Let $u\in \tilde{B}$, then there exists an element $S_7\in \mathfrak{H}\setminus \mathfrak{H}_t$ tangent to $\bar{B}\cup \tilde{B}$ in $u$. 
\end{prop}
\proof
Let now $u\in \tilde{B}$ and let $L_u=\langle \tilde{p},u\rangle$ be a line through $u$ such that $\tilde{p}\in X'$, $\tilde{p}\not=\tilde{t}$ and $L_u$ is tangent to $\tilde{B}$ in $u$.\\
Because of Lemma 3 there exists an element $Y\in \mathcal{S}$, $Y\subseteq \langle X,X'\rangle$, with $Y\not=X,X'$ such that $\tilde{p}\in \langle p,Y\rangle$ (i.e. $\{\tilde{p}\}= \left \langle p,Y \right\rangle\cap X'$) and $\left \langle p,Y \right\rangle\cap \langle t,r\rangle=\emptyset.$\\
Since $\langle p,Y\rangle\cap S_3=\langle p,Y\rangle\cap \langle r,\tilde{q},t\rangle$ and $\langle p,Y\rangle\cap \langle r,t \rangle=\emptyset$, the intersection $\langle p,Y\rangle\cap \langle r,\tilde{q},t\rangle$ is just a point, say $m_1$ and we have $m_1\not\in \langle t,r\rangle$.
Therefore there exists $m_2\in (L\setminus \{s\})\cap V$ such that the line $\langle m_1,m_2\rangle$ is imaginary, i.e. tangent to $\mathfrak{K}(r,V)$ at the point $m_2$ and hence has empty intersection with $\bar{B}$.
Since in $\Pi_3$ there exists a plane through $u,m_2,\{Y\} $, there exists $Z\in \mathcal{S}$ such that $u\in \langle m_2,Z,Y\rangle $.\\
Let us consider $S_7:=\langle p,m_2,Z,Y\rangle \in \mathfrak{H}$, since $Y\subseteq \langle X,X'\rangle$ and $X\not \subseteq\langle Z,Y\rangle$ we have that $X'\not \subseteq \langle Z,Y\rangle$ and hence $S_7\not \in \mathfrak{H}_t$.
Now, we want to prove that $S_7$ is tangent to $\bar{B}\cup \tilde{B}$ in $u$.
For Proposition 3 we have that $S_7\cap S_3$ is a line and hence, since $m_1,m_2 \in S_7$ we have: $S_7\cap S_3=\langle m_1,m_2\rangle$. Therefore:
$$S_7\cap \bar{B}=(S_7\cap S_3) \cap \bar{B}=\langle m_1,m_2\rangle \cap \bar{B}=\emptyset. $$
On the other hand, since $S_7\not \in \mathfrak{H}_t$, for Lemma 5, $S_7\cap \langle \tilde{B}\rangle$ is a line not contained in $\Sigma$ and hence, since $u, \tilde{p}\in S_7$, we have $S_7\cap \langle \tilde{B}\rangle=L_u$. Therefore:
$$S_7 \cap \tilde{B}=(S_7\cap \langle \tilde{B}\rangle)\cap \tilde{B}=L_u\cap \tilde{B}=\{u\}.$$
Hence we can find a tangent element of $\mathfrak{H}$ to all $u\in\tilde{B}$.
\endproof
Summing up the previous results we get:
\begin{teo}
$\bar{B}\cup \tilde{B} $ is a minimal and non planar $\mathfrak{H}$-blocking set. Therefore the cone $B:=\mathfrak{K}(p,\bar{B}\cup \tilde{B})$ is a minimal blocking set of $\Pi_3=PG(3,q^6)$.
\end{teo}
We note that, being non planar it is impossible to obtain this example using the classical MP construction.\\
Now we see that the blocking set $B$ appears not to be in the MPS class.
\begin{teo}
$B=\mathfrak{K}(p,\bar{B}\cup\tilde{B})$ is a minimal blocking sets of $\Pi_3=PG(3,q^6)$ not equivalent to any blocking sets of class MPS.
\end{teo}
\proof
First of all, let us evaluate the cardinality of $B=\mathfrak{K}(p,\bar{B}\cup \tilde{B})$.\\
We have that $\bar{B}=\mathfrak{K}(r,V\setminus L\cup \{s\})$ and hence: 
$$|\bar{B}\setminus \Sigma|=q^2(|V\setminus L\cup \{s\}|-|V\cap \Sigma|)=q^2((q^2+1)-1)=q^4.$$
Since $\tilde{B}$ is the union of three non parallel affine planes that share the same point we have that:
$$|\tilde{B}|=3q^4-3q^2+1.$$
Therefore we have that:
 $$|B|=|\mathfrak{K}(p,\bar{B}\cup \tilde{B})|=q^2(|(\bar{B}\cup \tilde{B})\setminus \Sigma|)+1=q^2(3q^4-3q^2+1+q^4)+1=4q^6-3q^4+q^2+1.$$
Let $q=p^e$ and let us suppose $B$ is obtained using the MPS construction starting from $PG(3n,q'=p^t) \mbox{ and } n\geq 2 $, hence $6e=nt$. Then we have seen (cf. [10], pag. 100) that $(p^t)^{n-1}|(|B|-1)$.\\
Since $|B|=q^2(4q^4-3q^2+1)+1$ and $p\not|4q^4-3q^2+1$, we have that $p^{t(n-1)}|q^2=p^{2e}$,
therefore we get $t(n-1)\leq 2e$ and $6e=t(n-1)+t\leq 2e+t $ which means that $4e\leq t=\frac{6e}{n}$. Hence we have the contradiction that $n=1$.
\endproof

\vspace{0.5cm}
\noindent
Dipartimento di Matematica e Fisica, Universit\`a Roma 3\\
Largo San Leonardo Murialdo 1,  00146,\ Roma, Italy\\
e-mail: {\bf costa@mat.uniroma3.it}


\begin{thebibliography}{10}
\bibitem{1} L. Bader, G. Lunardon, Desarguesian spreads,
\emph{Ricerche Mat.} 60 (2011), 15-37.
\bibitem{1Bis} A. Barlotti, J. Cofman, Finite Sperner Spaces Constructed from
Projective and Affine Spaces, \emph{Abh. Math. Semin. Univ. Hamb.} 40 (1974), 231-241.
\bibitem{2} S. G. Barwick, L. R. A. Casse, C. T. Quinn, The
Andr\'e/Bruck and Bose representation in $PG(2h,q)$: unitals and Baer
subplanes \emph{Bull. Belg. Math. Soc.} 7 (2000), 173-197.
\bibitem{3} L. M. Batten, \emph{Combinatorics of Finite Geometry},
Cambridge University Press, 2d ed., Cambridge, 1997.
\bibitem{4} A. Beutelspacher, F. Mazzocca, Blocking sets in infinite
projective and affine spaces, \emph{J. Geom.} 28, No 2 (1987),
112-116.
\bibitem{5} A. Blokhuis, P. Sziklai, T. Sz\H{o}ny, Blocking Sets in
Projective Spaces, \emph{Current Research Topics in Galois Geometry},
Nova Science Publishers Inc. (2011), 61-84.
\bibitem{5Tris} G. Lunardon, Linear k-blocking sets,
\emph{Combinatorica} 4 (2001), 571-581.
\bibitem{5bis} G. Lunardon, Normal spreads,
\emph{Geom. Dedicata} 75 (1999), 245-261.
\bibitem{6} G. Marino, O. Polverino, Ovoidal blocking sets and maximal
partial ovoids of Hermitian varieties, \emph{Des. Codes Cryptogr.} 56
(2010), 115-130.
\bibitem{7} F. Mazzocca, O. Polverino, L. Storme, Blocking Sets in
$PG(r, q^n)$, \emph{Des. Codes Cryptogr.} 44 (2007), 97-113.

\bibitem{8} F. Mazzocca, O. Polverino, Blocking sets in $PG(2, q^n)$
from cones of $PG(2n, q)$, \emph{J. Algebr. Comb.} 24 (2006), 61-81.
\bibitem{9} F. Mazzocca, G. Tallini, On the non existence of
blocking-sets in $PG(n,q)$ and $AG(n,q)$, for large enough $n$,
\emph{Simon Stevin} 1 (1985), 43-50.
\bibitem{10} F. Mazzocca, Geometrie di Galois, Codici, Disegni,
\emph{Appunti del corso di geometria combinatoria}, Napoli.
\bibitem{11} J. V. Neumann, O. Morgenstern, \emph{Theory of games and
economic behavior}, 2d ed.,Princeton University Press, Princeton 1947.
\bibitem{11Bis} O. Polverino, Linear sets in finite projective spaces,
\emph{Dis. Math.} 310 (2010), 3096-3107.
\bibitem{12} M. Richardson, On Finite Projective Games,
\emph{Proc. Amer. Math. Soc.} Vol. 7, No. 3 (June 1956), 458-465.
\bibitem{13} L. S. Shapley, A Value for $n$-person Games, \emph{In
H.W. Kuhn and A.W. Tucker (eds.), Contributions to the theory of games
II (AM-28)}, Princeton University Press, Princeton 1953.
\end{thebibliography}
\end{document}